\documentclass{amsart}
\usepackage{amssymb}

\def\xxx#1:{\medskip\noindent{\bf #1:}}
\def\dream#1:{\medskip\noindent{\shadowbox{\sf #1:}}}

\def\ldots{\hbox{$\mathinner {\ldotp \ldotp \ldotp }$}}

\newtheorem*{therevised}{The Revised GCH Theorem}
\renewcommand{\l}{\lambda}
\newcommand{\m}{\mu}
\renewcommand{\k}{\kappa}
\newcommand{\w}{\omega}
\renewcommand{\th}{\theta}

\newcommand{\pa}{\forall}

\newcommand{\pAa}{\rightarrow}

\newcommand{\pl}{\vert}

\newcommand{\ps}{\subseteq}

\newcommand{\pu}{\cup}

\newcommand{\px}{\times}
\newcommand{\py}{\exists}

\newcommand{\sP}{\mathcal{P}}

\newcommand{\ga}{\mathfrak{a}}

\newcommand{\gb}{\mathfrak{b}}

\newcommand{\gc}{\mathfrak{c}}

\newcommand{\xN}{\mathbb{N}}

\newcommand{\xQ}{\mathbb{Q}}

\newcommand{\xR}{\mathbb{R}}

\newcommand{\xZ}{\mathbb{Z}}

\newcommand{\haw}{\aleph_\omega}
\newcommand{\ha}[1]{\aleph_{#1}}

\newcommand{\hb}{\beth}
\newcommand{\lk}{\langle}
\newcommand{\rk}{\rangle}

\newtheorem{Thesis}{Thesis}

\newtheorem{definition}[Thesis]{Definition}

\newtheorem{theorem}[Thesis]{Theorem}

\title{You Can Enter Cantor's Paradise!}
\author{Saharon Shelah}

\begin{document}

\maketitle

\begin{center}
{\tt shelah@math.huji.ac.il}
\bigskip

\end{center}

I will try to use a spiralic presentation
returning to the same points on higher
levels hence repeating ourselves, so that a reader lost somewhere,
will not go away empty handed. Also I will assume essentially no
particular knowledge and I will say little on the history to which many great
mathematicians contributed.

\section{ Hilbert's first problem}
Recall (Cantor):
\begin{itemize}
\item We say that two sets $A,B$ are
{\em equinumerous} (or {\em equivalent\/}) if there is a one-to-one and
onto mapping from $A$ onto $B$;
\item The {\bf C}ontinuum {\bf H}ypothesis, CH, is the following
statement:

every infinite set of reals is either
equinumerous with the set ${\mathbb
Q}$ of rational numbers, or is
equinumerous with the set ${\mathbb R}$ of
all reals;
\item For a set $X$, let $\sP(X)$
denote its power set, i.e., the
set of all subsets of $X$.

The {\bf G}eneralized {\bf C}ontinuum {\bf H}ypothesis, GCH, is the
statement asserting that for every infinite set $X$, every subset $Y$ of
the power set ${\mathcal P}(X)$ is either equinumerous with a subset of $X$,
or is equinumerous with ${\mathcal P}(X)$ itself.
\end{itemize}

I think this problem is better understood
in the context of:
\subsection{Cardinal Arithmetic}

Recall (Cantor), that we call two sets $A,B$ equivalent
(or equinumerous)
if there
is a one-to-one mapping from $A$ onto $B$; the number of elements of
$A$ is the equivalence class of $A$ denoted by $\pl A\pl$, we call
it also {\em the power\/} or {\em the cardinality\/} of $A$. Having
defined infinite numbers, we can naturally ask ourselves what
is  the natural  meaning of  the arithmetical operations and   the
order.
There can be little doubt concerning  the order:

\begin{itemize}
\item $\pl A \pl \le \pl B\pl $ iff $A$ is equivalent to some subset of $B$.
\end{itemize}

We know that

\begin{itemize}
\item
any two  infinite cardinals are  comparable
so it is really a linear order

\item any cardinal $\l$ has a successor  $\l^+$, which means that

\item  $\l <\m \Leftrightarrow \l^+ \le \m$.
\end{itemize}
Well, but what about the arithmetical operations?
There are  natural definitions for the basic operations:
\begin{itemize}
\item addition (such that $\pl A\pu B\pl=\pl A\pl+\pl B\pl$ when $A,B$
are disjoint),
\item multiplication (such that $\pl A\px B\pl=\pl A\pl\px\pl B\pl$),
\item exponentiation (such that $\pl A\pl^{\pl B \pl}=\pl {}^B A \pl$,
where ${}^B A=\{f:f$ a function from $B$ to $A\}$).
\end{itemize}

A mathematician is allowed to choose his definitions and give
them ``nice" names,  but do those operations have any laws?
Are there interesting theorems about them?
For the second  questionit is a big  yes but is  outside
 the scope of this article. For the first question, the
answer is clear cut: all the  usual
\emph{equalities} hold, that is,  addition and multiplication
satisfy the commutative, associative, and distributive laws and
their infinite parallels. Also  for exponentiation, e.g.
 $(\l^\m)^\k = \l^{\m \px \k}, \l^\m \px \l^\k
= \l^{\m+\k}$. However this does not hold for
 the inequalities.  For every infinite cardinal $\l$ we
have $\l=
\l+1$. This should not surprise us; it is to be expected that allowing
inifinite numbers
 will ``cost" us some losses,
(as extending $\xN$,  the integers,  to rationals
``costs" us  the existence of
successor and the proof by induction; this is very clear in a posteriory
wisdom, of which we all have a lot).
Cantor was going around
asking: are there more points in the plane than on the line? People
answered him: don't you see that there are? But it is false, the line
and the plane are equinumerous.

In fact we can totally understand addition and multiplication, as the
following very nice rules holds for infinite numbers:
\begin{itemize}
\item $\mu + \l = \max(\mu, \l)$
\item $\mu \px \l = \max(\mu, \l)$
\end{itemize}
School children would have loved such arithmetic!


You may wonder: is this not too good?  Maybe all infinite numbers are equal
so this arithmetic is not so interesting?
But Cantor showed that $2^\l>\l$, meaning in particular that there are more
reals than natural numbers; noting that he called the number of natural
numbers $\ha0$, this is  the first infinite cardinal and showed that the
number of reals is $2^{\ha0}$.

Recall that every infinite number $\l$ has a successor, one bigger than
it  but smaller or equal to any bigger number, and it is denoted by $\l^+$.

Now mathematicians  tend to conjecture
that things are nice and well
understood. So, having only two natural
operations to increase a cardinal,
what is more natural than to conjecture
that those two operations, $\l^+$ and
$2^\l$ are equal.
Also mathematicians tend to conjecture either that whatever they cannot prove
may fail, or whatever they cannot build counterexample to is true;  and being
unable to construct an intermediate cardinal between $\l $ and $2^\l$
(e.g. $\ha0$ and $2^{\ha0}$) it is natural to conjecture that there is
nothing between them. This is

\subsection{Hilbert's first problem,  general version} The
``generalized continuum hypothesis'', or  GCH, says: for every
infinite number $\m$, its power
$2^\m$ is its successor $\m^+$.

The interest is that if GCH holds, then not only addition and
multiplication are easy, but also exponentiation is easy:
 for infinite cardinals $\l$, $\k$
(on $cf({\lambda})$ see below):

\[
\l^\k=
\begin{cases}
\l &  \text{ if } \k<cf(\l)\cr
\max(\ \k^+, \l^+)& \text{ if } \k\geq cf(\l)\cr
\end{cases}
\]

\textbf{Dream}: Find the laws of (infinite) cardinal exponentiation.

It has been assumed that if we understand cardinal arithmetic, that is
(taking for granted the understanding of addition and multiplication)
understand the behavior of exponentiation, we will generally  understand
set theory much better, and so solve problems from many branches of
mathematics in full generality.

\section{
Proven ignorance:  showing that we cannot know!}

\xxx The continuum problem:  How many real numbers are there?

\xxx Cantor proved:  There are {\em more}
reals than rationals.
(``There is no bijection from $\xR$ onto  the rationals $\xQ$'')

Recall that the continuum hypothesis (CH) says:  yes, more, but barely.
 Every set $A
\subseteq {\mathbb R}$ is either countable or
equinumerous with $\mathbb R$.

\xxx G\"odel proved:  Perhaps CH holds.

\xxx Cohen proved:   Perhaps CH does not hold.

\xxx G\"odel:  CH cannot be refuted.  Moreover, the {\em generalized
continuum hypothesis} may hold, in fact it holds if we restrict to
the  class $L$ of  constructible sets (\cite{godel}).

This universe $L$ satisfies all the axioms of set theory, and in addition
the generalized continuum hypothesis. The class $L$ can be described as the
minimal family of sets you absolutely must have as soon as you have all the
ordinals=order types of linear orders which are well ordered, i.e., every
non empty subset has a first element. We shall not deal with this here,
and do not touch on metamathematical matters in general.
\medskip

\xxx Cohen: You cannot prove that all sets are constructible, and you cannot
even prove the weaker statement ``CH'' (\cite{cohen}).

Cohen discovered the method of  {\em forcing\/},
and used it to prove this
``independence'' result; he ``fattened''  the universe of set theory; not
surprising in a posteriory wisdom. Again, this is not our topic.

Easton showed that there are no more rules
than the classical ones  if we
restrict ourselves to the so called regular
 cardinals (they include $\ha0$
and all successor cardinals; the
classical laws are: $2^\l>\l$ and
$cf(2^\l)>\l$; an infinite cardinal
$\l$ is regular if $cf (\l) = \l $,
on $cf$ see later).

Concerning the remaining cardinals, the so called singulars,
completing the theorem for them was thought of as a technical
problem.  The first such cardinal is $\ha\w=\sum\limits_{n=0,1,2,\ldots}
\aleph_n$, were $\ha0$ is the number of natural numbers and $\aleph_{n+1}$
is the successor of $\aleph_n$.

Very surprisingly, in the mid-seventies some rules were discovered by
Silver (\cite{silver}), and by Galvin and Hajnal (\cite{galvinhajnal}).
Let
$\aleph_{\omega_1}$ be the first
cardinal below which there are
 uncountably  (i.e.   $> \ha0$)  many cardinals. Call $\l$ strong
limit if $\m<\l\pAa 2^\m<\l$.
 If $\aleph_{\omega_1}$ is a strong limit
cardinal then below $2^{\aleph_{\omega_1}}$
there are at most $2^{\ha1}$
cardinals.

There were more works, but the general opinion was that what is left is just
showing fully we cannot prove anything more. E.g.  in '86  Leo Harrington told
me: \emph{``Cardinal arithmetic? Yes, it had been a great problem,
but now ..."} Not clear to me why, even if  the only thing left is proving
everthing is independent, this is not a major problem,
but this is irrelevant
here. However
 I would like to stress that  the
independence results help us to  discover new good theorems by
discarding many fruitless directions.

The book (\cite{pcf}) for which I am honoured by the Bolyai Prize
is based on:
\begin{Thesis}
[``Treasures are waiting for you''] There are
 many laws of (infinite)
cardinal arithmetic concerning exponentiation;
they look meagre as we
have concentrated on $2^\l$, but if we
 deal with relatively small
exponent and large base, then
there is much to be said.
\end{Thesis}

I think that though  GCH has not been
seriously considered as a true axiom, it has
influenced the way we thought about the problem,
so traditionally set theorists
concentrate on $2^\l$; but actually
it seem reasonable that on  the case
 $\l^\k$ with $\k \ll \l$
which is closer to finite products
(about which we know everything)
we will be able to say more.

{\bf We wonder}: What is the simplest open
 case of cardinal arithmetic?

Clearly we have to consider countable
products i.e.
 with the index set being the
set of natural numbers (as finite
products are easy),  or if you prefer, consider
 exponentiation of the form
$\l^{\aleph_0}$,

Now, $\aleph_0$ is the first infinite cardinal
(the cardinality of the set
${\mathbb Q}$ of rationals and the
 cardinality of  the set $\xN$
of natural numbers),
$2^{\aleph_0}=\aleph_0^{\aleph_0}$ is called  the
continuum, on which we know everything
(i.e., we know that we may not know
more). Moreover, let $\ha{1}$ be the
successor of $\ha{0}$, $\ha{2}$ be
the successor of $\ha{1}$ and generally
$\aleph_n$ be the $n$-th uncountable
cardinal (i.e $>\ha0$). Now it is not hard to prove that
$\aleph_n^{\aleph_0}=\max(\aleph_n,2^{\aleph_0})$.

So the first non trivial case is the (infinite) product of  those
numbers: $\prod\limits_n \aleph_n$, or equivalently,
\begin{equation}
\text{ what is  } \aleph_\omega^{\aleph_0}?
\end{equation}
where  $\aleph_\omega$ is the  sum of the $\ha n$'s.

If the continuum (i.e.  $2^{\ha0}$) is above all the $\aleph_n$, then this
product is equal to the continuum; so assume that the continuum is one of
them.

Let $\aleph_{\omega_n}$ be the first cardinal
 below which there are
$\aleph_n$ infinite cardinals.

\begin{theorem}
$\prod\limits_n \aleph_n < \aleph_{\omega_4}$ when the product is not
$2^{\ha0}$.
\end{theorem}


You may think this is a typographical error (in fact almost all who saw it
for  the first time were convinced this is a typographical error) and we
still do not know:

\textbf{Dream/Question}:  Why the hell is it four? Can we replace it by one?
Is 4 an artifact of the proof or  the best possible bound?

I think the four looks strange because we are looking at the problems from
a not so good perspective.

\section{pcf theory}
Close to my heart is

\begin{Thesis}
Cardinal arithmetic is loaded with
 independence results because we ask
the ``wrong" questions. The ``treasures''
thesis above is not enough; we
should replace cardinality by
 cofinality, a notion explained below (pcf theory).
More fully,  the unclarity comes from the interaction of two phenomena:
the values of $2^\l$ for $\l$ regular on which we know all the rules
(see above) and the ``cofinality arithmetic" where there is much  to be
said.
\end{Thesis}

As an illustration, let us look
again at $\haw^{\ha 0}$.

Look at the family $[\aleph_\w]^{\aleph_0}$
of countable subsets of
any set of cardinality $\aleph_\w$, e.g.
$\aleph_\w$ by
 the usual convention that
$\haw$ itself is such a set.
 It is partially ordered by inclusion.  Now,
instead of asking
about its cardinality as above
(we know that the number of countable
subsets of $\haw$ is
$\haw^{\ha 0}$),
 we {\bf ask} what is the
minimal number of members so that
 any other is included in at least one of
them,
and we call this number its cofinality, denoted by
$cf([\aleph_\omega]^{\aleph_0})$?
The theorem quoted above really says

\begin{theorem}
$cf([\aleph_\w]^{\aleph_0})<\aleph_{\omega_4}$.
\end{theorem}

This is meaningful even if the continuum, $2^{\ha 0}$ is large.
This exemplifies that even if one restricts oneself
to sets of reals only, set theory is not changed much
and in particular, cardinal arithmetic reinterprets as above
(even if we
restrict ourselves to simply defined sets of reals
(with arbitrary maps)).

We now pay some debts, defining $cf (\l)$ and regular cardinals.

\begin{definition}
\begin{enumerate}
\item A cardinal number
$\k$ is {\em regular},
if: whenever $A$ is of size $\k$,
 $A=\bigcup\limits_{i\in I} A_i$,
and all $A_i$ are of smaller cardinality than $A$, then $I$ must be at least
of size $\k$.
{\em (A set of size $\kappa$ cannot be written as a union of ``few''
``small'' sets.)}

\item Otherwise $\k$ is called singular;
 $cf(\k)$ is the size of the
smallest set $I$ that can appear in $(*)$ above.

\item For a partial order $P$ let $cf (P)$, its cofinality,
 be  the minimal cardinal
$\k$ such that   some subset $Q$ of $P$
of cardinality $\l$, $Q$ is cofinal in $P$, i.e.
$(\pa x\in P) (\py y\in Q ) [ x\le_P y ]$.
\end{enumerate}
\end{definition}

Note that successor cardinals, $\l^+$, are regular,
and the first singular
cardinal is $\aleph_\w$; regular limit cardinals
which are uncountable
 are ``large'',
and  the cofinality of any linear order with no
last elements is a regular cardinal.

So replacing $\l^\k$ by $cf([\l]^\k)$ we get
a much more ``robust'' theory,
were there are more answers and less ``we cannot answer''.

In fact, we are driven further. Suppose we consider linear orders $L_t$ for
$t\in T$ and assume that $cf (L_t)$ is bigger than the cardinality of  $T$.
We can look at the product, $\prod\limits_{t\in T} L_t$ ordered by:  $f\leq g$
iff $(\pa t\in T )(f(t)<_{L_t} g(t))$. This is not a linear order though it is
not grossly not so, because the number of factors is small compared to the
cofinality of the factors, as assumed. We can try to ``localize'' as done in
other  cases in mathematics; e.g. for analysing $\xZ$, we may like to
analyse what occurs when we have just one prime, so we consider  the
$p-$adics; sometime we may deduce information on $\xZ$ by looking at all those
 completions (and  the reals). Here we like to make  the order linear; so let
 let $I$ be a maximal ideal of the Boolean algebra of subsets of $T$, then
define {\bf linear} order $<_I$ on $\prod\limits_{t\in T} L_t$ by

$f<_I g$\quad iff\quad $\{t\in T:\neg [f(t)< g(t)]\}\in I$

So $cf (  \prod\limits_{t\in T} L_t , <_I )$ is a regular cardinal.

Now the product stops giving us a specific result and instead gives us a
spectrum:

\begin{definition} For a set $\ga=\{\l_t:t\in T\}$ of regular cardinals
each bigger than the cardinality of $T$, let $pcf(\ga)$ be the
set of cardinals of the form $cf(\prod\limits_{t\in T} L_t  , <_I )$,
where each $L_t$ is a linear order of cofinality $\l_t$ and $I$ is a maximal
ideal on  the Boolean algebra of subsets of $T$.
\end{definition}

This may remind you of what we do for rings for which decomposition to primes
does not behave as in the integers.
We have "the pcf theorem":

\begin{theorem} For $\ga=\{\l_t:t\in T\}$ as above:

(1)\quad $pcf(\ga)$ is a set of regular cardinals
 extending  $\ga$,
with at most $2^{|T|}$
   elements,

(2) \quad $pcf(\ga)$ has a largest element
denoted by $\max pcf(\ga)$.

(3)\quad The cofinality of the partial order
$\prod\limits_{t\in T}\l_t$ is
   $\max pcf(\ga)$.

(4)\quad For every $\th\in pcf(\ga)$ \mbox{there is a set
$\gb_\theta\ps\ga$} such that: for any maximal
 ideal $I$ on $T$ the
\mbox{cofinality of }
$(\prod\limits_{t\in T}
  \l_t ,<_I)$ is $\min\{\th\in pcf(\ga):\gb_\theta\notin I\}$.
\end{theorem}

Note that in part (1),  the number of maximal ideals on $T$
is $2^{2^{\pl T\pl}}$, so it gives some information (though not clear if
 the best possible one). Now this operation, $pcf$ has various rules, from
 them we can derive the bound $\aleph_{\omega_4}$ for
$cf([\aleph_\omega]^{\aleph_0},\subseteq)$.

Such rules are

{\bf Rule (local behavior):} If $\gb$ is a subset of $ pcf( \ga )$
(with both $\ga,\gb$ sets of regular cardinals
bigger than the number of members)
{\bf then}  any member $\l$ of $ pcf( \gb )$
belongs to $ pcf( \gc )$ for some subset $\gc$
of $\gb$ of cardinality smaller or equal to
the cardinality of $\ga$.

{\bf Rule (convexity):} $ pcf\{\ha n :n \ge 1 \}$
is an initial segment
of the set of successor cardinals.

{\bf Rule (continuity):} if
$L$ is a linear order,
 $\lk \l_i:i\in L \rk$ is
an increasing continuous sequence of cardinals bigger
than the cardinality of
$L$ and $\l= \Sigma_{i\in L} \l_i$ and $cf (\lambda)>\ha0$
then for some closed unbounded subset $C$ of $L$
we have $\max  pcf (\{ cf (\l_i ) : i\in C\}) = \l^+$.

Now we can explain better the inequality
 $cf ([\haw ]^{\ha0}, \subseteq )$:
  we look at  the set $pcf (\{ \aleph_n : n= 1,2,3,\dots \})$,
 on it we have a
linear order (the order on  the cardinals) and
 a closure operation
(which is pcf
 itself), this is topology of a special kind; we have enough
rules such that if  the size of  the set is too large then we get a
contradiction.


\section{ GCH Revisited}
We can interpret pcf theory as  a positive solution of Hilbert's
first problem, after considering  the inherent limitations.
We may interpret the problem more strictly ``cardinal arithmetic is
easy", and we present now a positive solution
from \cite{sh460}, discussed in details in its introduction, explaining why
reinterpreting GCH
under the known restrictions, we can prove that it holds almost always.
We define a variant of exponentiation,  which gives a different ``slicing"
of the GCH, specifically represent it as an equality  on
two cardinals ($\l^{\lk \k \rk} =\l$), and present a theorem saying that
the equality holds ``almost always".

\begin{definition} (1) For $\k<\l$ regular let $\l^{ \lk\k \rk}$,
 the `revised' power of
$\l$ by $\k$, be the minimal cardinality of a family
$\sP$ of subsets of $\l$ of cardinality $\k$
such that any other such set is the union of $<\k$  sets from $\sP$.

(2) For $\k<\l$ regular let $\l^{ [\k ]}$,
 the revised power of
$\l$ by $\k$, is the minimal cardinality of a family
$\sP$ of subsets of $\l$ of cardinality $\k$
such that any other such set is a subset of the union of $<\k$  sets from
$\sP$.
\end{definition}

{\bf Remark:} (1) GCH is equivalent to: for every regular $\l>\k$
we have $\l^{ \lk \k \rk} = \l$;

(2) Note that $\l^{[\k ]} \le  \l^{\lk \k \rk }  \le \l^{[\k ]}  +2^\k$
hence for $\l \ge 2^\k$ the two revised powers are equal, so below it
does not matter which version we use.

(3) So a weak version of GCH is: for ``most'' pairs  $(\l,\k)$
of regular cardinals we have $\l^{ [\k ]} = \l$.

{\bf Notation:} Let $ \beth_0=\ha0$,
$\hb_{n+1}= 2^{\hb_n}$, $\hb_\w = \sum_{n<{\omega}} \hb_n$.

\begin{therevised}
For any  $\l>\hb_\w$, for any $\k<\hb_\w$
large enough, we have $\l^{[\k ]} = \l$.
\end{therevised}


\begin{thebibliography}{9}

\bibitem{cohen} P.~J.~Cohen:
	{\it Set Theory and the Continuum Hypothesis},
	Benjamin, 1966.

\bibitem{galvinhajnal} F.~Galvin, A.~Hajnal:
	Inequalities for cardinal powers,
	{\it Annals of Mathematics}, {\bf 101}(1975), 491--498.

\bibitem{godel} K.~G\"odel:
	{\it The consistency of the axiom of choice and the generalized
        continuum-hypothesis with the axioms of set theory},
	Princeton University Press, 1940.

\bibitem{pcf} S.~Shelah: {\it Cardinal Arithmetic}, {Oxford University Press},
	1994

\bibitem{sh460} S.~Shelah: The Generalized Continuum Hypothesis revisited,
	{\it Israel Journal of Mathematics}, {\bf 116}(2000), 285--321.

\bibitem{silver} J.~Silver: On the singular cardinals problem,
	in: {\it Proceedings of the International Congress of Mathematics},
	Vancouver, 1974, volume I, 265--268.

\end{thebibliography}
\end{document}